\documentclass[a4paper, 12pt, reqno]{amsart}

\usepackage{amsfonts, amssymb, amsmath, eucal, verbatim, amsthm, amscd, enumerate}
\usepackage[all]{xy}

\newtheorem{thm}{Theorem}
\newtheorem*{thm*}{Theorem}
\newtheorem{lemma}[thm]{Lemma}
\newtheorem{prop}[thm]{Proposition}
\newtheorem*{prop*}{Proposition}
\newtheorem{cor}[thm]{Corollary}
\theoremstyle{definition}
\newtheorem{defn}[thm]{Definition}
\newtheorem{condition}[thm]{Condition}

\newcommand{\R}{\mathbb R}
\renewcommand{\H}{\mathbb H}
\newcommand{\Z}{\mathbb Z}
\newcommand{\Q}{\mathbb{Q}}
\newcommand{\D}{\mathbb D}
\newcommand{\sign}{\mathrm{sign}}
\renewcommand{\th}{^\mathrm{th}}
\newcommand{\NN}{\mathcal{N}}
\newcommand{\Or}{\mathrm{Or}}
\DeclareMathOperator{\rank}{rank}
\DeclareMathOperator{\interior}{int}

\newcommand{\BGMR}{MR1608567}
\newcommand{\BGM}{MR1293878}
\newcommand{\Plucker}{MR1288305}
\newcommand{\KSONE}{MR932303}
\newcommand{\KSTWO}{MR1094466}
\newcommand{\ZILLER}{ZILLERPREPRINT}
\newcommand{\CAYLEYTHEOREM}{MR1871828}
\newcommand{\THESIS}{HepworthThesis}

\author{Richard A. Hepworth}
\title{The topology of certain $3$-Sasakian $7$-manifolds.}
\address{School of Mathematics\\ University of Edinburgh}
\email{r.a.hepworth@ed.ac.uk}

\begin{document}

\begin{abstract}
We calculate the integer cohomology ring and stable tangent bundle of a family of compact, $3$-Sasakian $7$-manifolds constructed by Boyer, Galicki, Mann, and Rees.
Previously only the rational cohomology ring was known.
The most important part of the cohomology ring is a torsion group that we describe explicitly and whose order we compute.
There is a surprising connection with the combinatorics of trees.
\end{abstract}

\maketitle

\section{Introduction.}\label{Introduction_Section}

\subsection{Background and results.}
In \cite{\BGMR} Boyer, Galicki, Mann, and Rees constructed an infinite family of compact $3$-Sasakian $7$-manifolds by applying the technique of $3$-Sasakian reduction \cite{\BGM} to torus actions on spheres.
The rational cohomology ring of these manifolds was computed, showing that the second Betti number can be made arbitrarily large.
This in turn showed the existence of compact Einstein manifolds with positive scalar curvature and arbitrarily large total Betti number.

In this article we compute the integer cohomology of the manifolds constructed by Boyer, Galicki, Mann, and Rees.
We are able to immediately conclude that there are infinitely many homotopy types in each \emph{rational} homotopy type of the manifolds.
Also, for the manifolds $M$ in this family, the torsion cohomology group $H^4(M)$ is bounded below in terms of the second Betti number:
\[|H^4(M)|\geqslant(b_2(M)+2)^{b_2(M)}.\]

Our results represent the first stage of a homeomorphism and diffeomorphism classification of these manifolds, along the lines of the classifications \cite{\KSONE}, \cite{\KSTWO} made by Kreck and Stolz.
(The methods of Kreck and Stolz do not apply to the manifolds here, but the necessary generalisation of the Kreck-Stolz theory was made by the author \cite{\THESIS}.)
Such a classification would shed light on the extent to which $3$-Sasakian structures (and more generally, Einstein structures) restrict the topology of manifolds which possess them. 
In particular, it would be relevant to the question of how many different $3$-Sasakian structures can appear on the same manifold.
(Chinburg, Escher, and Ziller \cite{\ZILLER} have exhibited a single manifold with two non-isometric $3$-Sasakian structures.)\\[10 pt]

In order to give precise statements of our main results, recall that for each $(k+2)\times k$ integer matrix $\Omega$ satisfying Condition \ref{Reduction_Condition} below, Boyer, Galicki, Mann, and Rees \cite{\BGMR} constructed a compact 3-Sasakian manifold $S^7_\Omega$.
\begin{condition}\label{Reduction_Condition}
Let $\Delta_{pq}$, for distinct $p,q\in\{1,\ldots,k+2\}$,  denote the determinant of the matrix obtained by deleting the $p\th$ and $q\th$ rows of $\Omega$.
\begin{enumerate}
\item $\Delta_{pq}\neq 0$ for all $p\neq q$.
\item For each $p$, the $\Delta_{pq}$ for $q\neq p$ have greatest common divisor $1$.
\end{enumerate}
\end{condition}

The following two theorems summarise the existing results on the cohomology of the $S^7_\Omega$.
The first is obtained from \cite[Theorem E]{\BGM}, and the second is \cite[Theorem 4.2]{\BGMR}.

\begin{thm*}
In the case $k=1$, the manifold $S^7_\Omega$ has integer cohomology of the following form:
\[\begin{array}{c|cccccccc}
n & 0 & 1 & 2 & 3 & 4 & 5 & 6 & 7 \\
H^n(S^7_\Omega) & \Z & 0 & \Z & 0 & \Z/r & \Z^k & 0 & \Z
\end{array}\]
where $r=|\Delta_{12}||\Delta_{13}|+|\Delta_{21}||\Delta_{23}|+|\Delta_{31}||\Delta_{32}|$, and $H^4(S^7_\Omega)$ is generated by the square of a generator of $H^2(S^7_\Omega)$.
\end{thm*}

\begin{thm*}
For all $k$, $S^7_\Omega$ has rational cohomology of the following form:
\[\begin{array}{c|cccccccc}
n & 0 & 1 & 2 & 3 & 4 & 5 & 6 & 7 \\
H^n(S^7_\Omega;\Q) & \Q & 0 & \Q^k & 0 & 0 & \Q^k & 0 & \Q
\end{array}\]
\end{thm*}

Our main result, Theorem \ref{Cohomology_Thm}, extends both of these theorems.
Before stating it we establish some notation.
The space $S^7_\Omega$ is constructed as the base space of a $T^k$-bundle, and so we may define $x_i\in H^2(S^7_\Omega)$ to be the first Chern class associated to the $i\th$ factor of $T^k$.
Let $v_p=\sum_{i=1}^k\Omega_{pi}x_i\in H^2(S^7_\Omega)$, and for distinct $p,q$ let $\epsilon^p_q=(-1)^{p+q}\sign(p-q)\sign(\Delta_{pq})\in\{\pm 1\}$.

\begin{thm}\label{Cohomology_Thm}
The cohomology of $S^7_\Omega$ is as follows,
\[\begin{array}{c|cccccccc}
n & 0 & 1 & 2 & 3 & 4 & 5 & 6 & 7 \\
H^n(S^7_\Omega) & \Z & 0 & \Z^k & 0 & G_\Omega & \Z^k & 0 & \Z
\end{array}\]
where $H^2(S^7_\Omega)$ is freely generated by the $x_i$ and $G_\Omega$ is generated by the cup-products $x_ix_j$, subject only to the relations 
\[\epsilon^p_q v_pv_q+\epsilon^q_r v_qv_r+\epsilon^r_pv_rv_p=0\]
for distinct $p,q,r$.
This $G_\Omega$ is a torsion group of order
\begin{equation}\label{Order_Equation}\sum|\Delta_{t_1s_1}|\cdots|\Delta_{t_{k+1}s_{k+1}}|,\end{equation}
where the summand indexed by $t_1,s_1,\ldots,t_{k+1},s_{k+1}$ is included if and only if the graph on vertices $\{1,\ldots,k+2\}$ with edges $\{t_i,s_i\}$ is a \emph{tree}.
\end{thm}

Since we have assumed each $|\Delta_{pq}|$ to be nonzero, the sum \eqref{Order_Equation} is at least as large as the number of trees on the vertex set $\{1,\ldots,k+2\}$.
A theorem of Cayley \cite{\CAYLEYTHEOREM} states that the number of trees on a vertex set $A$ is precisely $|A|^{|A|-2}$.
We therefore have:
\begin{cor}
$|H^4(S^7_\Omega)|\geqslant(b_2(S^7_\Omega)+2)^{b_2(S^7_\Omega)}.$
\end{cor}

In \cite{\BGMR}, it was shown that a matrix $\Omega$ of the form
\[\Omega=\left(\begin{array}{ccc}
1   &        &     \\
    & \ddots &     \\
    &        & 1   \\
a_1 & \cdots & a_k \\
b_1 & \cdots & b_k 
\end{array}\right)\]
satisfies Condition \ref{Reduction_Condition} so long as:
\begin{enumerate}
\item all pairs $a_i$ $b_i$ are coprime;
\item one never has $a_i=\pm a_j$, $b_i=\pm b_j$ for distinct $i,j$.
\end{enumerate}
Matrices $\Omega$ of this form can be constructed with arbitrarily large $a_i$, $b_i$, for example by letting $a_1,\ldots,b_k$ be distinct primes.
Further, it is clear from Theorem \ref{Cohomology_Thm} that $H^4(S^7_\Omega)$ will have size at least as large as each of the $|a_i|$, $|b_j|$.
Therefore for each $k$, one can construct an infinite sequence of $\Omega$ for which the groups $H^4(S^7_\Omega)$ are of unbounded size.
We therefore obtain:
\begin{cor}
There are infinitely many homotopy types in each \emph{rational} homotopy type of the $S^7_\Omega$.
\end{cor}

Our second main result is a computation of the stable tangent bundle and hence of the Pontrjagin classes of the $S^7_\Omega$.
\begin{thm}\label{Pontrjagin_Thm}
\[TS^7_\Omega\oplus\epsilon_\R^{4k+1}\cong 2\bigoplus_{p=1}^{k+2}l_p\]
where, for $1\leqslant p\leqslant (k+2)$, $l_p$ is the complex line bundle over $S^7_\Omega$ with first Chern class $v_p$.
Consequently $$p_1(S^7_\Omega)=2\sum_{p=1}^{k+2}v_p^2.$$
\end{thm}

Note that there are some differences in notation between the present article and the article \cite{\BGMR} where the $S^7_\Omega$ were constructed.
The matrices $\Omega$ here are the transposes of the $\Omega$ used in \cite{\BGMR}.
The manifold we denote $S^7_\Omega$ is written in \cite{\BGMR} as $\mathcal{S}(\Omega^t)$.

In the next two parts of this introductory section we review the construction of the manifolds $S^7_\Omega$, and we outline the rest of this article.

\subsection{The construction of the $S^7_\Omega$.}\label{Construction_Section}
We will review the construction of the manifolds $S^7_\Omega$.
These were originally constructed in \cite{\BGMR} using the technique of \emph{$3$-Sasakian reduction} developed in \cite{\BGM}.

Let $k\geqslant 2$ be an integer.
We regard the sphere $S^{4k+7}$ as the unit vectors in $\H^{k+2}$.
There is a canonical $3$-Sasakian metric on $S^{4k+7}$, under which the $T^{k+2}$-action defined by
\[(\mathbf{t},\mathbf{x})\mapsto\mathbf{t}\mathbf{x}=(t_px_p)_{p=1}^{k+2},\ \ \ \mathbf{t}\in T^{k+2},\ \mathbf{x}\in S^{4k+7},\]
consists of $3$-Sasakian isometries.
(Here we regard $T^1$ as the unit complex numbers, and therefore as a subgroup of the unit quaternions $S^3$.)
The $3$-Sasakian moment map associated to this action of $T^{k+2}$ is
\begin{gather}
\mu\colon S^{4k+7}\to \Im\H^{k+2},\label{mu_Equation}\\
\mu(\mathbf{x})=(\bar x_pix_p)_{p=1}^{k+2},\nonumber
\end{gather}
where $\Im\H$ denotes the imaginary quaternions.

Any $(k+2)\times k$ integer matrix $\Omega$ defines a homomorphism $T^k\to T^{k+2}$.
By combining this homomorphism with the action of $T^{k+2}$ on $S^{4k+7}$, we obtain an action of $T^k$ on $S^{4k+7}$ which consists of $3$-Sasakian isometries.
Explicitly, the action is given by
\[(\mathbf{t},\mathbf{x})\mapsto\mathbf{t}\mathbf{x}=(t_1^{\Omega_{p1}}\cdots t_{k+2}^{\Omega_{p\,{k+2}}}x_p)_{p=1}^{k+2},\]
and the associated $3$-Sasakian moment-map is
\begin{gather}
\mu_\Omega\colon S^{4k+7}\to\Im\H^k,\label{muOmega_Equation}\\ 
\mu_\Omega(\mathbf{x})= ({\textstyle\sum_{p=1}^{k+2}}\bar x_pix_p\Omega_{pj})_{j=1}^k.\nonumber
\end{gather}
That is, $\mu_\Omega(\mathbf{x})=\Omega^t\mu(\mathbf{x})$.

\begin{thm*}[{\cite[2.14]{\BGMR}}]
Suppose that Condition \ref{Reduction_Condition} holds for the matrix $\Omega$.
Then $3$-Sasakian reduction can be applied to the action of $T^k$ on $S^{4k+7}$, yielding a compact, $3$-Sasakian $7$-manifold
\[S^7_\Omega=S^{4k+7}/\negthinspace/T^k.\]
That is to say that $\mathbf{0}$ is a regular value of $\mu_\Omega$, $T^k$ acts freely on the submanifold $\mu_\Omega^{-1}(\mathbf{0})$ of $S^{4k+7}$, and the resulting quotient
\[S^7_\Omega=T^k\backslash \mu_\Omega^{-1}(\mathbf{0})\]
is a compact $3$-Sasakian manifold of dimension $7$.
\end{thm*}

\subsection{Outline of the article.}
We will now outline the methods by which Theorem \ref{Cohomology_Thm} and Theorem \ref{Pontrjagin_Thm} will be proved.

Recall the \emph{Borel Construction}.
Let $X$ be a topological space with a left-action by a topological group $G$.
The {Borel Construction} of $G$ on $X$ is the space $EG\times_GX=G\backslash(EG\times X)$.
Throughout this article we will write
\begin{align*}
\pi\colon & EG\times_G X\to BG,\\
Q\colon & EG\times_G X\to G\backslash X,
\end{align*}
for the maps induced by projection onto the first and second factors of $EG\times X$.
The map $\pi$ is the projection in a bundle with fibres $X$.
When $X$ is a principal $G$-bundle, $Q$ is a homotopy-equivalence.

For simplicity we write
\[N_\Omega=\mu_\Omega^{-1}(\mathbf{0}).\]
Since $T^k$ acts freely on $N_\Omega$, $Q\colon ET^k\times_{T^k}N_\Omega\to T^k\backslash N_\Omega=S^7_\Omega$ is a homotopy-equivalence.
Therefore the Serre spectral sequence of the bundle
\begin{equation}\label{BorelOne_Equation}
N_\Omega\to ET^k\times_{T^k}N_\Omega\xrightarrow{\pi}BT^k
\end{equation}
can be used to calculate $H^\ast(S^7_\Omega)$, \emph{if} we have information on $H^\ast(N_\Omega)$ and the spectral sequence differentials.
This information is obtained as follows.

Since $\mu_\Omega$ is invariant with respect to the $T^{k+2}$-action on $S^{4k+7}$, $T^{k+2}$ acts on $N_\Omega$.
The $T^{k+2}$ moment-map has strong properties, and using these we are able to understand the orbit structure and singular sets of the $T^{k+2}$-space $N_\Omega$.
With this information we compute $H^\ast(ET^{k+2}\times_{T^{k+2}}N_\Omega)$, which in turn is used in the Serre spectral sequence of 
\begin{equation}\label{BorelTwo_Equation}
N_\Omega\to ET^{k+2}\times_{T^{k+2}}N_\Omega\xrightarrow{\pi} BT^{k+2}
\end{equation}
to partially compute $H^\ast(N_\Omega)$ and the spectral sequence differentials that affect it.

Finally, the map $\Omega\colon T^k\to T^{k+2}$ is used to compare the Serre spectral sequences of \eqref{BorelOne_Equation} and \eqref{BorelTwo_Equation}.
The information provided by this comparison is sufficient to compute $H^\ast(S^7_\Omega)$.

The computation of the stable normal bundle of $S^7_\Omega$ is less involved.
Since $N_\Omega$ is the $\mathbf{0}$-set of the $T^{k+2}$-invariant map $\mu_\Omega\colon S^{4k+7}\to\Im\H^{k}$, its normal bundle in $S^{4k+7}$ is trivial with trivial $T^{k+2}$-action.
The $T^{k+2}$-equivariant stable normal bundle of $S^{4k+7}$ itself is easily understood, and so we can compute the $T^{k+2}$-equivariant normal bundle of $N_\Omega$.
Reducing with respect to $T^k$, we obtain the stable normal bundle of $S^7_\Omega$.\\[10 pt]

The rest of the article is organised as follows:
In Section \ref{Action_Section} we investigate $N_\Omega$ as a $T^{k+2}$-space, and in particular we give a complete description of the quotient space and the singular orbits.
In Section \ref{Borel_Section} we use the results of Section \ref{Action_Section} to compute $H^\ast(ET^{k+2}\times_{T^{k+2}} N_\Omega)$ in terms of $\pi\colon ET^{k+2}\times_{T^{k+2}}N_\Omega\to BT^{k+2}$.
The proof of one lemma that is entirely algebraic and rather different in nature from the rest of Section \ref{Borel_Section} is deferred to Section \ref{Cycle_Section}.
In Section \ref{GOmega_Section} \emph{part} of Theorem \ref{Cohomology_Thm}, essentially the result on the order of $G_\Omega$, is proved as a separate proposition; it is a lengthy exercise in algebra and the theory of trees.
Finally, in Sections \ref{Cohomology_Section} and \ref{Pontrjagin_Section} we prove Theorem \ref{Cohomology_Thm} and Theorem \ref{Pontrjagin_Thm} respectively.\\[5 pt]

\noindent\textbf{Acknowledgements:}
Thanks are due to my supervisor Professor Elmer Rees, for his advice and support during the writing of my thesis, from which this article is derived.
Thanks also to the Engineering and Physical Sciences Research Council of the United Kingdom, which supported me during my doctoral studies.

\section{The action of $T^{k+2}$ on $N_\Omega$.}\label{Action_Section}
From now on we fix an integer $k\geqslant 2$ and a $(k+2)\times k$ matrix $\Omega$ satisfying Condition \ref{Reduction_Condition}.
Recall that $N_\Omega=\mu^{-1}_\Omega(\mathbf{0})$ admits a $T^{k+2}$-action.
The object of this section is to understand the orbit space and singular set of this action.
The crucial results, which will be used in Section \ref{Borel_Section}, are Corollary \ref{NOmegaQuotient_Cor}, Proposition \ref{Orbits_Prop}, and Proposition \ref{NpOmega_Prop}.

We begin with the following lemma, which shows that the moment-map $\mu\colon S^{4k+7}\to \Im\H^{k+2}$ fully describes the $T^{k+2}$-orbits on $S^{4k+7}$.

\begin{lemma}\label{mu_Lemma}
Recall that $\mu\colon S^{4k+7}\to\Im\H^{k+2}$ is given by $\mu(\mathbf{x})=(\bar x_pix_p)_{p=1}^{k+2}$.
Give $\Im\H^{k+2}$ the norm $\|\mathbf{x}\|=\sum_{p=1}^{k+2}|x_p|$, and use the norm to form the unit sphere $S^{3k+5}\subset\Im\H^{k+2}$.
Then:
\begin{enumerate}
\item $\mu(S^{4k+7})=S^{3k+5}.$
\item The fibres $\mu^{-1}(\mathbf{s})$, $\mathbf{s}\in S^{3k+5}$, are precisely the $T^{k+2}$-orbits in $S^{4k+7}$.
\end{enumerate}
\end{lemma}
\begin{proof}
Let $\phi\colon\H\to\Im\H$ be the map defined by $\phi(u)=\bar uiu$.
One can verify the following:
\begin{enumerate}
\item $\phi$ is surjective.
\item The fibres of $\phi$ are the $T^1$-orbits $T^1u$.
\item That $|\phi(u)|=|u|^2$.
\end{enumerate}
Since $\mu$ is obtained by applying $\phi$ in each coordinate, the lemma follows immediately from the three properties of $\phi$.
\end{proof}

Since $\mu_\Omega=\Omega^t\circ\mu$ and $N_\Omega=\mu^{-1}_\Omega(\mathbf{0})$, $\mu|N_\Omega$ takes values in $\ker(\Omega^t\colon\Im\H^{k+2}\to\Im\H^k)$, which is a copy of $\Im\H^2\cong\R^6$, and inherits a norm from $\Im\H^{k+2}$.
Combining this with Lemma \ref{mu_Lemma} we obtain:

\begin{cor}\label{NOmegaQuotient_Cor}
The moment-map $\mu\colon S^{4k+7}\to\Im\H^{k+2}$ induces a homeomorphism 
\[T^{k+2}\backslash N_\Omega\to S^{3k+5}\cap\ker(\Omega^t\colon\Im\H^{k+2}\to\Im\H^k)=S^5.\]
\end{cor}

Now we establish some notation and use it describe the structure of the singular orbits in $N_\Omega$.

\begin{defn}\label{Singular_Defn}\hfill
\begin{enumerate}
\item For $1\leqslant p\leqslant (k+2)$ we denote by $N_\Omega^p\subset N_\Omega$ the subset of elements whose $p\th$ entry vanishes. 
\item $N^\bullet_\Omega=\bigcup_{p=1}^{k+2} N^p_\Omega$.
\item We denote by $T^1_p\subset T^{k+2}$ the subgroup of elements whose entries are all $1$ except for the $p\th$.
\item We write
\[T^1_{\epsilon^p}=\{(\mathbf{t},s)\in T^{k+2}\times S^3|s\in T^1,\ t_p=1,\ t_q=s^{\epsilon^p_q}\ \mathrm{for}\ q\neq p\}.\]
\end{enumerate}
\end{defn}

\begin{prop}\label{Orbits_Prop}
The torus $T^{k+2}$ acts freely on $N_\Omega-N^\bullet_\Omega$, and points in $N_\Omega^p$ have stabiliser $T^1_p$.
\end{prop}
\begin{proof}
An element $\mathbf{t}\in T^{k+2}$ fixes $\mathbf{x}\in N_\Omega$ if and only if $t_p=1$ whenever $x_p\neq 0$.
Since the elements of $N_\Omega-N^\bullet_\Omega$ have no entry equal to $0$, $T^{k+2}$ acts freely on $N_\Omega-N^\bullet_\Omega$.
Also, it was shown in the proof of \cite[Lemma 2.15]{\BGMR} that no element of $N_\Omega$ has more than one entry equal to $0$, so that points in $N^p_\Omega$ have stabiliser precisely $T^1_p$.
\end{proof}

The group $S^3$ of unit quaternions acts on $S^{4k+7}\subset\H^{k+2}$ by right-multiplication, and respects $\mu_\Omega$ in the sense that $\mu_\Omega(\mathbf{x}s)=s^{-1}\mu_\Omega(\mathbf{x})s$.
It follows that there is a left-action of the group $T^{k+2}\times S^3$, defined by $(\mathbf{t},s)\mathbf{x}=\mathbf{t}\mathbf{x}s^{-1}$, on $N_\Omega$ and on each $N^p_\Omega$.

\begin{prop}\label{NpOmega_Prop}\hfill
\begin{enumerate}
\item $N^p_\Omega$ is a submanifold of $N_\Omega$.
Its normal bundle is the trivial bundle $N^p_\Omega\times\H$, and $T^{k+2}$ acts on this normal bundle by $\mathbf{t}(\mathbf{x},h)=(\mathbf{t}\mathbf{x},t_ph)$.
\item As a $T^{k+2}\times S^3$ space, each $N^p_\Omega$ has the form $(T^{k+2}\times S^3)/(T^1_pT^1_{\epsilon^p})$; consequently $T^{k+2}\backslash N^p_\Omega\approx S^3/T^1=S^2$.
\end{enumerate}
\end{prop}

\begin{proof}
Let $\pi_p\colon N_\Omega\to\H$ be the map $\mathbf{x}\mapsto x_p$.
This is smooth, and $\pi_p(\mathbf{t}\mathbf{x})=t_p\pi_p(\mathbf{x})$, $\pi_p^{-1}({0})=N^p_\Omega$.
Therefore to prove the first part we need only show that $0$ is a regular value of $\pi_p$.

One can easily calculate $D_\mathbf{v}\mu_\Omega\colon T_\mathbf{v}S^{4k+7}\to T_{\mu_\Omega(\mathbf{v})}$ and verify that, for $\mathbf{n}\in N^p_\Omega$,
\[T_\mathbf{n}N_\Omega=\ker D_\mathbf{n}\mu_\Omega\subset T_\mathbf{n}S^{4k+7}\subset T_\mathbf{n}\H^{k+2}=\H^{k+2}\]
contains the subspace of vectors in $\H^{k+2}$ whose components all vanish except for the $p\th$.
This subspace is just a copy of $\H$, and $D_\mathbf{n}\pi_p\colon T_\mathbf{n}N_\Omega\to T_0\H=\H$ maps it isomorphically onto $\H$.
Therefore $D_\mathbf{n}\pi_p$ is surjective, so that $0$ is a regular value of $\pi_p$.
This proves the first part of the lemma.

To prove the second part of the lemma we must describe $N^p_\Omega$ as a space under the $T^{k+2}\times S^3$ action.
To do this, we will construct an explicit element $\mathbf{v}^p$ of $N^p_\Omega$ and then compute its orbit and stabiliser.

Define $\mathbf{v}^p\in S^{4k+7}$ as follows:
Set $e^p_q=1$ if $\epsilon^p_q=1$, and $e^p_q=j$ if $\epsilon^p_q=-1$.
Fix $p$ and form the vector in $\H^{k+2}$ whose $p\th$ entry is $0$ and whose $q\th$ entry for $q\neq p$ is $e^p_q\sqrt{|\Delta_{pq}|}$.
Scale this vector by a positive real number to obtain $\mathbf{v}^p\in S^{4k+7}$.

Note that $\bar e^p_qie^p_q=\epsilon^p_qi$, so that $\mu(\mathbf{v}^p)\in\Im\H^{k+2}$ is a positive real multiple of the vector whose $p\th$ entry vanishes and whose $q\th$ entry for $q\neq p$ is $\epsilon^p_q|\Delta_{pq}|i=\Delta_{pq}(-1)^{p+q}\sign(p-q)i$.
Consequently, the $q\th$ entry of $\Omega^t\mu(\mathbf{v}^p)$ is a multiple of
\begin{multline*}
\sum_{r=1}^{p-1}(-1)^{p+r}\Delta_{pr}\Omega_{rq}-\sum_{r=p+1}^{k+2}(-1)^{p+r}\Delta_{pr}\Omega_{rq}\\
=(-1)^{p+k+2}\left[\sum_{r=1}^{p-1}(-1)^{k+2+r}\Delta_{pr}\Omega_{rq}+\sum_{r=p+1}^{k+2}(-1)^{r-1+k+2}\Delta_{pr}\Omega_{rq}\right].
\end{multline*}
The quantity in square brackets is the determinant of the $(k+1)\times(k+1)$ matrix obtained by appending to $\Omega$ a copy of its $q\th$ column and then deleting the $p\th$ row; this matrix has a repeated column, and so the above quantity vanishes.
Thus $\mu_\Omega(\mathbf{v}^p)=\Omega^t\mu(\mathbf{v}^p)=0$, so that $\mathbf{v}^p\in N_\Omega$.
Since $v^p_p=0$, we have $\mathbf{v}^p\in N^p_\Omega$ as required.

Suppose that $\mathbf{n}\in N^p_\Omega$.
We will show that $\mathbf{n}\in(T^{k+2}\times S^3)\mathbf{v}^p$.
Fix $q\neq p$.
Since $v^p_q$ and $n_q$ are both nonzero, we can find $r\in(0,\infty)$ and $s\in S^3$ such that $rs\mu(\mathbf{n})s^{-1}$ and $\mu(\mathbf{v}^p)$ have $q\th$ components equal.
Now $rs\mu(\mathbf{n})s^{-1}-\mu(\mathbf{v}^p)$ is an element of $\ker(\Omega^t\colon \Im\H^{k+2}\to \Im\H^k)$ with two components zero; since all $k\times k$ submatrices of $\Omega^t$ are nondegenerate it follows that $rs\mu(\mathbf{n})s^{-1}-\mu(\mathbf{v}^p)=0$.
Thus $r=1$ and $\mu(\mathbf{n}s^{-1})=s\mu(\mathbf{n})s^{-1}=\mu(\mathbf{v}^p)$, so that by Corollary \ref{NOmegaQuotient_Cor} the vectors $\mathbf{n}s^{-1}$ and $\mathbf{v}^p$ lie in the same $T^{k+2}$ orbit.
That is, $\mathbf{n}$ and $\mathbf{v}^p$ lie in the same $T^{k+2}\times S^3$ orbit, as claimed.

Finally we compute the $(T^{k+2}\times S^3)$-stabiliser of $\mathbf{v}^p$.
Since $v^p_p=0$, $T^1_p$ fixes $\mathbf{v}^p$.
Recall that $v^p_q$ is an element of $\R$ or $j\R$ according to whether $\epsilon^p_q$ is $1$ or $-1$ respectively, so that $v^p_qz=z^{\epsilon^p_q}v^p_q$; it follows that $T^1_{\epsilon^p}$ fixes $\mathbf{v}^p$.
Now suppose that $(\mathbf{t},s)\in T^{k+2}\times S^3$ fixes $\mathbf{v}^p$.
Necessarily $s\in T^1$, and so we can multiply $(\mathbf{t},s)$ by an element of $T^1_pT^1_{\epsilon^p}$ to obtain an element $(\mathbf{t}',1)$ in the stabiliser of $\mathbf{v}^p$ with $t'_p=1$.
It follows immediately that $\mathbf{t}'=1$ and so $(\mathbf{t},s)\in T^1_pT^1_{\epsilon^p}$.
\end{proof}

\section{The Borel construction of $T^{k+2}$ on $N_\Omega$.}\label{Borel_Section}
In this section we will prove the following theorem, which describes the cohomology of the Borel construction $ET^{k+2}\times_{T^{k+2}}N_\Omega$ in terms of the projection-map $\pi\colon ET^{k+2}\times_{T^{k+2}}N_\Omega\to BT^{k+2}$.
The proof is based entirely on the results of Section \ref{Action_Section}.

We will use the notation $[n]=\{1,\ldots,n\}$ for natural numbers $n$.
We write $H^\ast(BT^{k+2})=\Z[u_1,\ldots,u_{k+2}]$, where $u_p$ is the first Chern class associated to the $p\th$ factor $BT^1$ of $BT^{k+2}$, and set
\begin{equation}\label{rho_Equation}
\rho_{pqr}=\epsilon^p_qu_pu_q+\epsilon^q_ru_qu_r+\epsilon^r_pu_ru_p,
\end{equation}
for distinct $p,q,r\in[k+2]$.

\begin{thm}\label{BorelCohomology_Thm}
The $\rho_{pqr}$ lie in the kernel of
\[\pi^\ast\colon H^\ast(BT^{k+2})\to H^\ast(ET^{k+2}\times_{T^{k+2}}N_\Omega),\]
and $\pi^\ast$ induces an isomorphism
\[\tilde\pi^\ast\colon H^\ast(BT^{k+2})/\langle\rho_{pqr}\mid p,q,r\mathrm{\ distinct}\rangle\to H^\ast(ET^{k+2}\times_{T^{k+2}}N_\Omega).\]
\end{thm}

We now fix some notation.
Let $T^1$ act on $S^3$ in the usual way, so that $T^1\backslash S^3=S^2$.
We fix the generator $c_1\in H^2(BT^1)$, and we choose the generator $r$ of $H^2(S^2)$ which is compatible in the sense that $Q^\ast r=\pi^\ast c_1$.
That is, $r$ is the first Chern class of the circle-bundle $S^3\to S^2$.

\begin{defn} Recall from Proposition \ref{Orbits_Prop} and Proposition \ref{NpOmega_Prop} that the singular orbits of $T^{k+2}$ in $N_\Omega$ are disjoint and have union $N_\Omega^\bullet=\bigcup_{p=1}^{k+2}N^p_\Omega$, where each $N^p_\Omega$ is a submanifold of $N_\Omega$.
We define the associated Borel constructions:
\begin{align*}
\NN_\Omega&=ET^{k+2}\times_{T^{k+2}}N_\Omega,\\
\NN_\Omega^p&=ET^{k+2}\times_{T^{k+2}}N^p_\Omega,\\
\NN_\Omega^\bullet&=ET^{k+2}\times_{T^{k+2}}N^\bullet_\Omega={\textstyle\bigcup_{p=1}^{k+2}}\NN^p_\Omega.
\end{align*}
\end{defn}

The object of Theorem \ref{BorelCohomology_Thm} is to compute $H^\ast(\NN_\Omega)$.
The next three lemmas, whose proofs are given at the end of this section, compute $H^\ast(\NN^\bullet_\Omega)$, $H^\ast(\NN_\Omega,\NN^\bullet_\Omega)$, and the coboundary $\delta^\ast\colon H^\ast(\NN^\bullet_\Omega)\to H^\ast(\NN_\Omega,\NN^\bullet_\Omega)$.
The proof of Theorem \ref{BorelCohomology_Thm} then consists of combining the three computations above using the long exact sequence of the pair $(\NN_\Omega,\NN_\Omega^\bullet)$.
This relies on the entirely algebraic Lemma \ref{Ideal_Lemma}, whose proof is given in Section \ref{Cycle_Section}.

\begin{lemma}\label{NNDotCohomology_Lemma}
Recall from Proposition \ref{NpOmega_Prop} that $T^{k+2}\backslash N_p^\bullet=\bigcup_{p=1}^{k+2}S^2$; we write $r_p\in H^2(T^{k+2}\backslash N^\bullet_\Omega)$ for the usual generator in the $p\th$ summand.
Then
\begin{equation*}H^\ast(\NN^\bullet_\Omega)=\bigoplus_{p=1}^{k+2}\Z[U_p,Q^\ast r_p]/Q^\ast r_p^2,\label{NNOmegaDotCohomology_Equation}\end{equation*}
where $\pi^\ast u_p=U_p+\sum_{q\neq p}\epsilon^q_pQ^\ast r_q$.
With this description, the restriction $H^\ast(\NN^\bullet_\Omega)\to H^\ast(\NN^p_\Omega)$ is just projection onto the $p\th$ summand.
\end{lemma}

\begin{lemma}\label{PairIsomorphism_Lemma}
Recall from Corollary \ref{NOmegaQuotient_Cor} and Proposition \ref{NpOmega_Prop} that $T^{k+2}\backslash N_\Omega\approx S^5$ and $T^{k+2}\backslash N^p_\Omega\approx S^2$.
The map $Q^\ast\colon H^\ast(S^5,\bigcup_{p=1}^{k+2}S^2)\to H^\ast(\NN_\Omega,\NN^\bullet_\Omega)$ is an isomorphism.
In particular, $H^n(\NN_\Omega,\NN^\bullet_\Omega)=0$ unless $n=1,3,5$.
\end{lemma}

\begin{lemma}\label{Delta_Lemma}
The coboundary map $\delta^\ast\colon H^{\ast-1}(\NN_\Omega^\bullet)\to H^\ast(\NN_\Omega,\NN_\Omega^\bullet)$ is surjective.
Its kernel is described below.
\begin{enumerate}
\item $\delta^\ast\colon H^0(\NN_\Omega^\bullet)\to H^1(\NN_\Omega,\NN_\Omega^\bullet)$ has kernel $\sum 1_p=Q^\ast\sum 1_p$.
\item $\delta^\ast\colon H^2(\NN_\Omega^\bullet)\to H^3(\NN_\Omega,\NN_\Omega^\bullet)$ has kernel spanned by the $\pi^\ast u_p=U_p+\sum_{q\neq p}\epsilon^q_pQ^\ast r_q$.
\item $\delta^\ast\colon H^4(\NN_\Omega^\bullet)\to H^5(\NN_\Omega,\NN_\Omega^p)$ has kernel generated by the $U_p^2=\pi^\ast u_p^2$ and $U_pQ^\ast r_p-U_q Q^\ast r_q=\epsilon^p_q\pi^\ast u_pu_q$.
\item In all other degrees, $\ker(\delta^\ast)=H^\ast(\NN_\Omega^\bullet)$.
\end{enumerate}
\end{lemma}

Recall from \eqref{rho_Equation} that $\rho_{pqr}=\epsilon^p_qu_pu_q+\epsilon^q_ru_qu_r+\epsilon^r_pu_ru_p$ for $p,q,r$ distinct.
In the next section we will prove:
\begin{lemma}\label{Ideal_Lemma}
The ideal $\langle\rho_{p,q,r}\mid p,q,r\ \mathrm{distinct}\rangle$ contains all terms
\begin{gather*}
u_pu_qu_r,\\
\epsilon^p_qu_p^2u_q-\epsilon^p_ru^2_pu_r,\\
u_p^2u_q^2,
\end{gather*}
for $p,q,r$ distinct.
\end{lemma}

\begin{proof}[Proof of Theorem \ref{BorelCohomology_Thm}]
By Lemma \ref{Delta_Lemma}, $\delta^\ast\colon H^\ast(\NN_\Omega^\bullet)\to H^\ast(\NN_\Omega,\NN_\Omega^\bullet)$ is surjective, so that $H^\ast(\NN_\Omega)\to\ker(\delta^\ast)\subset H^\ast(\NN_\Omega^\bullet)$ is an isomorphism.

It is immediate that $H^0(\NN_\Omega)=\Z$ and that $H^{2i-1}(\NN_\Omega)=0$ for $i\geqslant 1$, so that the map $\tilde\pi^\ast$ of the statement is an isomorphism in degree $0$ and in odd degrees.
By Lemma \ref{NNDotCohomology_Lemma} and the second part of Lemma \ref{Delta_Lemma}, $H^2(\NN_\Omega)$ is free on the classes $\pi^\ast u_p$, so that $\tilde\pi^\ast$ is an isomorphism in degree $2$.
By Lemma \ref{NNDotCohomology_Lemma} and the third part of Lemma \ref{Delta_Lemma}, $H^4(\NN_\Omega)$ is spanned by the classes $\pi^\ast u_pu_q$ for $p,q\in[k+2]$, subject only to the relations 
\[\epsilon^p_q\pi^\ast u_pu_q+\epsilon^q_r\pi^\ast u_qu_r+\epsilon^r_p\pi^\ast u_ru_p=0,\]
so that $\tilde\pi^\ast$ is an isomorphism in degree $4$.

It remains to prove that $\tilde\pi^\ast$ is an isomorphism in degrees $2n$ for $n>2$.
In such degrees $H^{2n}(\NN_\Omega)\to H^{2n}(\NN_\Omega^\bullet)$ is an isomorphism and so, by Lemma \ref{NNDotCohomology_Lemma}, $H^{2n}(\NN_\Omega)$ is freely generated by the classes 
\[\pi^\ast u_1^n,\ldots,\pi^\ast u_{k+2}^n,\pi^\ast u_1^{n-1}u_{l_1},\ldots,\pi^\ast u_{k+2}^{n-1}u_{l_{k+2}},\]
where each $l_p$ is any choice of element in $[k+2]-\{p\}$.
However, by Lemma \ref{Ideal_Lemma}, the elements
\[u_1^n,\ldots,u_{k+2}^n,u_1^{n-1}u_{l_1},\ldots, u_{k+2}^{n-1}u_{l_{k+2}},\]
span $H^\ast(BT^{k+2})/\langle \rho_{pqr}\rangle$ in degree $2n$.
Therefore in degree $2n$ the map $\tilde\pi^\ast$ sends a spanning set to a basis, and must therefore be an isomorphism.
This completes the proof.
\end{proof}

\begin{proof}[Proof of Lemma \ref{NNDotCohomology_Lemma}]
Recall from Proposition \ref{NpOmega_Prop} that, as a $T^{k+2}$-space, $N^p_\Omega=(T^{k+2}\times S^3)/T^1_{\epsilon^p}T^1_p$, where 
\[T^1_{\epsilon^p}=\{(\mathbf{t},s)\in T^{k+2}\times S^3|s\in T^1,\ t_p=1,\ t_q=s^{\epsilon^p_q}\ \mathrm{for}\ q\neq p\}.\]
By writing $T^{k+1}_p\subset T^{k+2}$ for the subgroup of elements with $p\th$ entry $1$, and splitting $T^{k+2}$ as $T^1_p\times T^{k+1}_p$, we can simplify the description of $N^p_\Omega$ to $(T^{k+1}_p\times S^3)/T^1_{\epsilon^p}$.
Consequently
\begin{eqnarray}
\NN^p_\Omega
&=&ET^{k+2}\times_{T^{k+2}}(T^{k+1}_p\times S^3)/T^1_{\epsilon^p}\nonumber\\
&=&BT^1_p\times(ET^{k+1}_p\times_{T^{k+1}_p}(T^{k+1}_p\times S^3)/T^1_{\epsilon^p})\nonumber\\
&\cong&BT^1_p\times T^1_{\epsilon^p}\backslash(ET^{k+1}_p\times S^3),
\end{eqnarray}
where $T^1_{\epsilon^p}$ acts on $ET^{k+1}_p\times S^3$ via its inclusion in $T^{k+1}_p\times S^3$.
The final identification above was made using the identification
\[(ET^{k+1}_p\times_{T^{k+1}_p}(T^{k+1}_p\times S^3)/T^1_{\epsilon^p})\cong T^1_{\epsilon^p}\backslash(ET^{k+1}_p\times S^3)\]
where $[e,[\mathbf{t},s]]$ corresponds to $[\mathbf{t}^{-1}e,s^{-1}]$.
Recall from Proposition \ref{NpOmega_Prop} that $T^{k+2}\backslash N^p_\Omega=S^2$.
Under the above identification of $\NN^p_\Omega$,
\[Q\colon \NN^p_\Omega\to T^{k+2}\backslash N^p_\Omega=S^2\]
becomes the map
\[BT^1_p\times T^1_{\epsilon^p}\backslash(ET^{k+1}_p\times S^3)\to T^1\backslash S^3=S^2.\]
This map, when restricted to $b\times T^1_{\epsilon^p}\backslash(ET^{k+1}_p\times S^3)$ for any $b\in BT^1_p$, is a homotopy-equivalence, and so:
\begin{equation}\label{NNpCohomologyOne_Equation}H^\ast(\NN^p_\Omega)=\Z[U_p,Q^\ast r_p]/Q^\ast r_p^2,\end{equation}
where $U_p=\pi^\ast u_p$ and $r_p\in H^2(S^2)$ is just a copy of the canonical generator $r$.

From the commutative diagram
\[\xymatrix{
T^1_{\epsilon^p}\backslash(ET^{k+1}_p\times S^3)\ar[r] & BT^{k+1}_p\\
ET^1\times_{T^1}S^3 \ar[r]\ar[u]& BT^1,\ar[u]
}\]
where the vertical maps are induced by the group homomorphism $T^1\cong T^1_{\epsilon^p}\hookrightarrow T^{k+1}_p\times S^3$, $z\mapsto ((z^{\epsilon^p_q})_{q\neq p},z)$, it follows that
\begin{equation}\label{NNpCohomologyTwo_Equation}\pi^\ast u_q=\epsilon^p_q Q^\ast r_p\ \ \mathrm{for}\ q\neq p.\end{equation}
Combining \eqref{NNpCohomologyOne_Equation} and \eqref{NNpCohomologyTwo_Equation} for all $p$, we obtain the result.
\end{proof}

\begin{proof}[Proof of Lemma \ref{PairIsomorphism_Lemma}]
By Proposition \ref{NpOmega_Prop}, $N^p_\Omega$ is a submanifold of $N_\Omega$, and has a normal bundle preserved by the $T^{k+2}$-action.
Consequently $N^\bullet_\Omega$ is a strong equivariant deformation retract of some neigbourhood $M$ which is preserved by the $T^{k+2}$-action.
In turn, $\NN^\bullet_\Omega$ is a strong deformation retract of the neighbourhood $\mathcal{M}=ET^{k+2}\times_{T^{k+2}}M$ and $T^{k+2}\backslash N_\Omega^\bullet=\bigcup^{k+2}_{p=1}S^2$ is a strong deformation retract of the neighbourhood $T^{k+2}\backslash M$.
Therefore we have a commutative diagram:
\[\xymatrix{
H^\ast(\NN_\Omega,\NN^\bullet_\Omega)  &  H^\ast(\NN_\Omega,\mathcal{M})\ar[l]_{\cong}\ar[r]^\cong  &  H^\ast(\NN_\Omega-\NN^\bullet_\Omega,\mathcal{M}-\NN_\Omega^\bullet)\\
H^\ast(S^5,\bigcup S^2)\ar[u]^{Q^\ast}   &  H^\ast(S^5,T^{k+2}\backslash M)\ar[l]_\cong\ar[r]^-\cong\ar[u]^{Q^\ast}   &  H^\ast(S^5-\bigcup S^2,T^{k+2}\backslash M-\bigcup S^2)\ar[u]^{Q^\ast}
}\]
The horizontal maps on the left are isomorphisms by homotopy-invariance; the horizontal maps on the right are isomorphisms by excision.

We claim that the vertical map on the right is an isomorphism.
Using the diagram above, the result will follow from this claim.
Since $T^{k+2}$ acts freely on $N_\Omega-N^\bullet_\Omega$ and $M-N^\bullet_\Omega$, the maps $N_\Omega-N^\bullet_\Omega\to S^5-\bigcup S^2$ and $M-N^\bullet_\Omega\to T^{k+2}\backslash M-\bigcup S^2$ are principal $T^{k+2}$-bundles.
It follows, as in the discussion of the Borel construction at the start of this section, that the associated maps
\begin{align*}
Q\colon&\NN_\Omega-\NN^\bullet_\Omega=ET^{k+2}\times_{T^{k+2}}(N_\Omega-N^\bullet_\Omega)\to S^5-\bigcup S^2,\\
Q\colon&\mathcal{M}-\NN^\bullet_\Omega=ET^{k+2}\times_{T^{k+2}}(M-N^\bullet_\Omega)\to T^{k+2}\backslash M-\bigcup S^2,
\end{align*}
are homotopy-equivalences and so induce isomorphisms on cohomology.
The claim that the vertical map at the right of the commutative diagram above is an isomorphism now follows by the five-lemma, and this completes the proof.
\end{proof}

\begin{proof}[Proof of Lemma \ref{Delta_Lemma}.]
By Lemma \ref{PairIsomorphism_Lemma}, $H^n(\NN_\Omega,\NN^\bullet_\Omega)$ vanishes except when $n=1,3,5$.
Therefore to prove the lemma we need only show that, when $n=1,3,5$, the coboundary $\delta^\ast\colon H^{n-1}(\NN_\Omega^\bullet)\to H^n(\NN_\Omega,\NN^\bullet_\Omega)$ is surjective with kernel as described.

In the case $n=1$, the claim is immediate since $Q^\ast\colon H^0(\bigcup S^2)\to H^0(\NN_\Omega^\bullet)$ and $Q^\ast\colon H^1(S^5,\bigcup S^2)\to H^1(\NN_\Omega,\NN^\bullet_\Omega)$ are isomorphisms.

In the case $n=3$, the map $\delta^\ast$ is surjective because $H^3(\NN_\Omega,\NN_\Omega^\bullet)$ is free on the $Q^\ast\delta^\ast r_p=\delta^\ast(Q^\ast r_p)$.
The elements described do lie in $\ker(\delta^\ast)$ since they are the restrictions of elements of $H^2(\NN_\Omega)$; that these elements span $\ker(\delta^\ast)$ follows by using Lemma \ref{NNDotCohomology_Lemma} to compare ranks.

In the case $n=5$, the elements described do lie in $\ker(\delta^\ast)$, since they are the restrictions of elements in $H^4(\NN_\Omega)$.
Note that, by Lemma \ref{PairIsomorphism_Lemma}, $H^5(\NN_\Omega,\NN_\Omega^\bullet)$ is free on a single generator $Q^\ast t$, where $t\in H^5(S^5,\bigcup S^2)$ restricts to a generator of $H^5(S^5)$.
We claim that $\delta^\ast U_pQ^\ast r_p=\pm t$ for each $p$.
From this, $\delta^\ast$ is surjective, and the claim for $n=5$ will follow by using Lemma \ref{NNDotCohomology_Lemma} to compare ranks.

We will now prove the claim that $\delta^\ast(U_pQ^\ast r_p)=\pm t$ for each $p$.

Using the exponential map, which is $T^{k+2}$-equivariant, we consider $\bigcup_p N^p_\Omega\times\H$ to be embedded in $N_\Omega$ as a tubular neighbourhood of $N_\Omega^\bullet$.
For brevity we define
\begin{align*}
\bar N_\Omega&=\bigcup N^p_\Omega\times\D^4,& \bar \NN_\Omega&=ET^{k+2}\times_{T^{k+2}}\bar N_\Omega,\\
\widetilde N_\Omega&=\bigcup N^p_\Omega\times\interior\D^4,&   \widetilde\NN_\Omega&=ET^{k+2}\times_{T^{k+2}}\widetilde N_\Omega,\\
N_\Omega^\circ&=\bigcup N^p_\Omega\times S^3,&    \NN_\Omega^\circ&=ET^{k+2}\times_{T^{k+2}}N_\Omega^\circ.
\end{align*} 

The space $N_\Omega-\widetilde N_\Omega$ is obtained from $N_\Omega$ by deleting the open unit disc bundle from each of the embedded normal bundles $N^p_\Omega\times\H$.
In each of the resulting fibres $\H-\interior\D^4$ of $\bigcup N^p_\Omega\times\H-\tilde N_\Omega$, we can collapse the boundary sphere $S^3$ to a point and identify the quotient space $(H-\interior\D^4)/S^3$ with $\H$.
In this way we obtain an equivariant map $C\colon(N_\Omega-\widetilde N_\Omega,N_\Omega^\circ)\to(N_\Omega,N_\Omega^\bullet)$ which induces an isomorphism on cohomology.
Also, since $C$ is equivariant it induces a map on the associated quotient spaces and Borel constructions; we also denote this induced map by $C$.

Note that $T^{k+2}\backslash N_\Omega^\circ=\bigcup_{p=1}^{k+2}S^2\times S^2$, where in each component the first factor is $T^{k+2}\backslash N^p_\Omega$, and the second factor is $T^1_p\backslash S^3=S^2$.
Therefore, the generators in the $p\th$ summand $H^2(S^2\times S^2)$ of $H^2(\bigcup S^2\times S^2)$ are $C^\ast r_p$ and $s_p$, where $s_p$ is a copy of the canonical generator $r$ associated to the second factor $S^2$.
Note that $Q^\ast s_p=i_p^\ast\pi^\ast u_p$, where $i_p^\ast\colon H^\ast(\NN^p_\Omega)\to H^\ast(\NN^\bullet_\Omega)$ denotes the inclusion of the $p\th$ summand.

Now we have:
\begin{eqnarray}
C^\ast\delta^\ast(U_pQ^\ast r_p)
&=&\delta^\ast((C^\ast U_p)(Q^\ast C^\ast r_p))\nonumber\\
&=&\delta^\ast((i_p^\ast\pi^\ast u_p)(Q^\ast C^\ast r_p))\nonumber\\
&=&\delta^\ast((Q^\ast s_p)(Q^\ast C^\ast r_p))\nonumber\\
&=&Q^\ast \delta^\ast(s_p C^\ast r_p)\label{Nasty_Equation}
\end{eqnarray}
The element $s_pC^\ast r_p\in H^4(\bigcup S^2\times S^2)$ is the dual to the fundamental class of the $p\th$ component $S^2\times S^2$, and therefore $\delta^\ast(s_p C^\ast r_p)$ is dual to the fundamental class of the manifold-boundary pair $(S^5-\bigcup S^2\times\interior\D^3,\bigcup S^2\times S^2)$.
Since $C\colon(S^5-\bigcup S^2\times\int\D^3,\bigcup S^2\times S^2)\to(S^5,\bigcup S^2)$ induces an isomorphism on cohomology, this description of $\delta^\ast(s_p C^\ast r_p)$ means that $\delta^\ast(s_p C^\ast r_p)=\pm C^\ast t$.
Now by \eqref{Nasty_Equation}, $C^\ast\delta^\ast(U_pQ^\ast r_p)= \pm C^\ast Q^\ast t$ and, again since $C^\ast$ is an isomorphism, the claim follows.
\end{proof}

\section{The cycle.}\label{Cycle_Section}

The object of this section is to prove Lemma \ref{Ideal_Lemma}.
Recall from Section \ref{Introduction_Section} that we have defined signs
\[\epsilon^p_q=(-1)^{p+q}\sign(p-q)\sign(\Delta_{pq})\in\{\pm 1\}\]
for all distinct $p,q\in[k+2]$.

\begin{defn}\label{Adjacency_Defn}
We say that distinct $p,q\in[k+2]$ are \emph{adjacent} if $\epsilon^p_r\epsilon^q_r$ is independent of $r\in[k+2]-\{p,q\}$.
\end{defn}

\begin{lemma}\label{CycleExists_Lemma}
The graph on vertex set $[k+2]$ whose edges correspond to adjacent pairs is a cycle of length $(k+2)$.
\end{lemma}
\begin{proof}
Consider the vector $\mathbf{w}^p\in\R^{k+2}$ whose $p\th$ entry is $0$ and whose $q\th$ entry for $q\neq p$ is $\epsilon^p_q|\Delta_{pq}|$.
The vector $\mathbf{w}^p$ was shown in the proof of Proposition \ref{NpOmega_Prop} to be in the kernel of $\Omega^t$.
Moreover, since each $k\times k$ submatrix of $\Omega^t$ is nondegenerate, $\mathbf{w}^p$ actually spans the space $\{\mathbf{w}\in\ker\Omega^t\mid w_p=0\}$.

Since $\epsilon^p_q$ is the sign of the $q\th$ entry of $\mathbf{w}^p$, $p$ and $q$ are adjacent if and only if $\mathbf{w}^p$ and one of $\pm\mathbf{w}^q$ share the same sign in all but their $p\th$ and $q\th$ entries.
This is if and only if one of the line segments $t\mathbf{w}^p+(1-t)\mathbf{w}^q$, $t\mathbf{w}^p-(1-t)\mathbf{w}^q$, $t\in[0,1]$, does not meet any of the lines $\{\mathbf{w}\in\ker\Omega^t\mid w_r=0\}=\R\mathbf{w}^r$ for $r\neq p,q$.
Sketching $\ker\Omega^t\cong\R^2$ and the $(k+2)$ lines $\R\mathbf{w}^r$ in it, one sees that $p$ and $q$ are adjacent if and only if the two lines $\R\mathbf{w}^p$ and $\R\mathbf{w}^q$ are adjacent in the sketch.
The result now follows.
\end{proof}

\begin{defn}
We define $\Or_{pqr}=\epsilon^p_q\epsilon^q_r\epsilon^r_p$ for distinct $p,q,r\in[k+2]$.
\end{defn}

\begin{lemma}\label{Orientation_Lemma}
$\Or_{pqr}+\Or_{prs}+\Or_{psq}=\Or_{qrs}$ for all distinct $p,q,r,s\in[k+2]$.
\end{lemma}
\begin{proof}
First note that $\Or_{abc}$ is antisymmetric in $a$, $b$, $c$, and that consequently the claim holds for fixed $p,q,r,s$ if and only if it holds for each permutation of $p,q,r,s$.
Second note that, by the definition of adjacency, $\Or_{pqr}=\Or_{p'qr}$ whenever $p$ and $p'$ are adjacent.
As a consequence of these two facts, and of what we know from Lemma \ref{CycleExists_Lemma} about the structure of the graph that describes adjacency, we can without loss assume that the pairs $(p,q)$, $(q,r)$, and $(r,s)$ are all adjacent.
From this, we have $\Or_{pqr}=\Or_{pqs}=-\Or_{psq}$, and that $\Or_{prs}=\Or_{qrs}$, so that the required identity holds.
\end{proof}

\begin{proof}[Proof of Lemma \ref{Ideal_Lemma}]
One can immediately check, by expanding the $\rho_{abc}$ in terms of $u_a$, $u_b$, $u_c$ and using Lemma \ref{Orientation_Lemma}, that
\[\Or_{pqr}[\epsilon^q_s\epsilon^r_s u_p\rho_{qrs}+\epsilon^r_s\epsilon^p_su_q\rho_{rps}+\epsilon^p_s\epsilon^q_su_r\rho_{pqs}]=-u_pu_qu_r\]
for all distinct $p,q,r,s\in[k+2]$.
This shows that $u_pu_qu_r\in\langle \rho_{abc}\mid a,b,c\mathrm{\ distinct}\rangle$ for each $p,q,r$ distinct.
It now follows, by adding an appropriate multiple of $u_pu_qu_r$ to $u_p\rho_{pqr}$, that $\epsilon^p_qu_p^2u_q-\epsilon^p_ru^2_pu_r\in\langle \rho_{abc}\mid a,b,c\mathrm{\ distinct}\rangle$.
Finally, multiplying $\epsilon^p_qu_p^2u_q-\epsilon^p_ru^2_pu_r$ by $u_q$ and using the first result, one sees that $u_p^2u_q^2\in\langle \rho_{abc}\mid a,b,c\mathrm{\ distinct}\rangle$.
\end{proof}

\section{The group $G_\Omega$.}\label{GOmega_Section}
In this section we will prove part of Theorem \ref{Cohomology_Thm}, rewritten as Proposition \ref{GOmega_Proposition}.
The proposition is an essential step in the proof of Theorem \ref{Cohomology_Thm}, but is proved by means entirely different from the rest of that theorem.

\begin{defn}
For a finite set $A$, $T_A$ will denote the collection of \emph{trees} (that is, connected acyclic graphs) on the vertex set $A$.
By regarding an edge in a graph as an element of the set $A^{(2)}$ of unordered pairs of distinct elements in $A$, we can regard $T_A$ as a subset of the set $(A^{(2)})^{(|A|-1)}$ of unordered $(|A|-1)$-tuples of distinct elements in $A^{(2)}$.
\end{defn}

\begin{defn}\label{newxivp_Defn}
Let $x_i\in H^2(BT^k)$ denote the first Chern class associated to the $i\th$ factor $BT^1$ of $BT^k$, so that $H^\ast(BT^k)=\Z[x_1,\ldots,x_k]$.
For $p\in[k+2]$ set $v_p=\sum_i\Omega_{pi}x_i=B\Omega^\ast u_p$.
\end{defn}

\begin{prop}\label{GOmega_Proposition}
The group $G_\Omega=H^4(BT^k)/\mathrm{span}(\epsilon^p_qv_pv_q+\epsilon^q_rv_qv_r+\epsilon^r_pv_rv_p)$ is finite and has order
\[\sum_{\mathbf{t}\in T_{[k+2]}}\prod_{\{p,q\}\in\mathbf{t}}|\Delta_{pq}|.\]
\end{prop}

The symbols in Definition \ref{newxivp_Defn} were used in Section \ref{Introduction_Section} to denote elements of $H^2(S^7_\Omega)$, and $G_\Omega$ was defined in the statement of Theorem \ref{Cohomology_Thm} in terms of these elements of $H^2(S^7_\Omega)$.
It is immediate to check that $x_i,v_p\in H^2(S^7_\Omega)$ are the pullbacks under the classifying map $S^7_\Omega\to BT^k$ of $x_i,v_p\in H^2(S^7_\Omega)$, and that this identifies the two descriptions of $G_\Omega$.

\begin{lemma}\label{GCD_Lemma}
$\gcd_{\mathbf{t}\in T_{[k+2]}}({\textstyle \prod_{\{p,q\}\in\mathbf{t}}}|\Delta_{pq}|)=1$.
\end{lemma}
\begin{proof}
Suppose that the lemma were false.
Then there is a prime $P$ which divides ${\textstyle \prod_{\{p,q\}\in\mathbf{t}}}|\Delta_{pq}|$ for each $\mathbf{t}\in T_{[k+2]}$.

We first prove the following:
If distinct $p,q,r\in [k+2]$ are such that $P\mid\Delta_{pq},\Delta_{qr}$, then $P\mid\Delta_{pr}$.
For, by the second part of Condition \ref{Reduction_Condition} there is $s\in [k+2]-\{q\}$ such that $P\nmid \Delta_{qs}$.
The minor determinants of $\Omega$ satisfy the Pl\"ucker relations, which for our purposes take the following form:
For distinct $p_1,p_2,p_3,p_4\in[k+2]$,
\[\Delta_{p_1p_3}\Delta_{p_2p_4}\pm\Delta_{p_2p_3}\Delta_{p_1p_4}\pm\Delta_{p_1p_2}\Delta_{p_3p_4}=0.\]
(Adapted from \cite[Chapter VII, Section 5]{\Plucker}.)
It now follows that since $P$ divides $\Delta_{pq}\Delta_{rs}$ and $\Delta_{ps}\Delta_{qr}$, $P$ must also divide $\Delta_{qs}\Delta_{pr}$, and so $P$ divides $\Delta_{pr}$.

Since $P\mid\Delta_{pq}$ for at least one $\{p,q\}\in [k+2]^{(2)}$, we can find nonempty sets $C\subset [k+2]$ such that $P\mid\Delta_{pq}$ for all $\{p,q\}\in C^{(2)}$.
Now let $C$ be maximal amongst sets with this property.  (That is, if we enlarge $C$ the property fails.)
Obviously $C\neq [k+2]$, and $C\neq\emptyset$, we may choose $\mathbf{t}\in T_{[k+2]}$ such that no element of $\mathbf{t}$ lies in $C^{(2)}$ or $([k+2]-C)^{(2)}$.
Since $P\mid(\prod_{\{p,q\}\in\mathbf{t}}|\Delta_{pq}|)$, there exists $\{p,q\}\in\mathbf{t}$ such that $P\mid\Delta_{pq}$.
Without loss assume that $p\in [k+2]-C$ and $q\in C$.
For any $q'\in C-\{q\}$ we have, by our assumption on $C$, that $P\mid\Delta_{qq'}$, and since also $P\mid\Delta_{pq}$, we have that $P\mid\Delta_{pq'}$ by the property established at the beginning of the proof.
Thus $C\cup\{p\}$ has the same property as $C$ did, contradicting the maximality of $C$.
Thus the lemma is proved.
\end{proof}

\begin{lemma}\label{NewDescription_Lemma}
The group $H^4(BT^k)$ is spanned by the classes
\begin{equation*}\label{Vpq_Equation}
V_{pq}=\epsilon^p_qv_pv_q
\end{equation*}
for $p\neq q$.
The $V_{pq}$ satisfy the relations
\begin{gather}
V_{pq}=-V_{qp}\ \ \mathrm{for\ each\ }p,q,\label{VpqRelationOne_Equation}\\
\sum_{q\colon q\neq p}|\Delta_{pq}|V_{pq}=0\ \ \mathrm{for\ each\ }p,\label{VpqRelationTwo_Equation}
\end{gather}
and \emph{all} relations among the $V_{pq}$ are generated by \eqref{VpqRelationOne_Equation} and \eqref{VpqRelationTwo_Equation}.

That is, $H^4(BT^k)$ is isomorphic to the quotient $V/R$, where $V$ is the free Abelian group spanned by elements $V_{pq}=-V_{qp}$ for distinct $p,q\in[k+2]$, and $R$ is the Abelian group spanned by the
\begin{equation*}
R_p=\sum_{q\colon q\neq p}|\Delta_{pq}|V_{pq},
\end{equation*}
for $p\in[k+2]$
\end{lemma}

\begin{proof}[Proof of Lemma \ref{NewDescription_Lemma}]
We begin by showing that the $V_{pq}$ do span $H^4(BT^k)$, and that they do satisfy relations \eqref{VpqRelationOne_Equation} and \eqref{VpqRelationTwo_Equation}.
By Condition \ref{Reduction_Condition}, any $(k+1)$ rows of the matrix $\Omega$ span $\Z^k$, so that any $(k+1)$ of the $v_p=\sum_i\Omega_{pi}x_i$ span $H^2(BT^k)$.
Therefore the $V_{pq}=\epsilon^p_qv_pv_q$, for fixed $p$ and $q\neq p$ varying, span all multiples of $v_p$ in $H^4(BT^k)$.
Thus the $V_{pq}$, as $p$ and $q$ vary, span all multiples of \emph{all} $v_p$ and, again since the $v_p$ span $H^2(BT^k)$, the $V_{pq}$ span $H^4(BT^k)$.
Relation \eqref{VpqRelationOne_Equation} is immediate since $\epsilon^p_q=-\epsilon^q_p$.
Also, as in the proof of Lemma \ref{NpOmega_Prop}, the sum $\sum_{q\colon q\neq p}\epsilon^p_q|\Delta_{pq}|\Omega_{qi}$ vanishes for each $i$ and $p$, so that
\begin{eqnarray*}
\sum_{q\colon q\neq p}|\Delta_{pq}|V_{pq}
&=& v_p\sum_{q\colon q\neq p}\epsilon^p_q|\Delta_{pq}|v_q\\
&=& v_p\sum_{q\colon q\neq p}\epsilon^p_q|\Delta_{pq}|\sum_{i=1}^k\Omega_{qi}x_i\\
&=& v_p\sum_{i=1}^kx_i\sum_{q\colon q\neq p}\epsilon^p_q|\Delta_{pq}|\Omega_{qi}\\
&=&0,
\end{eqnarray*}
for each $p$, which proves relation \eqref{VpqRelationTwo_Equation}.

We now show that the group described in the second paragraph of the statement is indeed isomorphic to $H^4(BT^k)$, which will complete the proof.
Note that $\sum_p R_p=0$ since $|\Delta_{pq}|V_{pq}=-|\Delta_{qp}|V_{qp}$ for each $p,q$.
Therefore $\rank (R)\leqslant(k+1)$, but in fact $\rank (R)=(k+1)$ since $R_p$ for $p\neq 1$ contains $V_{1p}$ with nonzero coefficient and no other $V_{1q}$ terms.
Therefore $V/R$ has rank $\binom{k+2}{2}-(k+1)=k+\binom k2=\rank H^4(BT^k)$, so that $V/R\to H^4(BT^k)$ is a surjection between Abelian groups of equal rank.
Since $H^4(BT^k)$ is free, it will suffice to show that $V/R$ is also free.
We will do this by showing that the torsion subgroup $\mathrm{torsion}(V/R)$ of $V/R$ vanishes.

For each $\mathbf{t}\in T_{[k+2]}$ consider the splitting $V=V_\mathbf{t}\oplus V'_\mathbf{t}$, where $V_\mathbf{t}$ is the subgroup spanned by the $V_{pq}$ for $\{p,q\}\in\mathbf{t}$, and $V'_\mathbf{t}$ is the subgroup spanned by the $V_{pq}$ for $\{p,q\}\not\in \mathbf{t}$.
Write $\pi_\mathbf{t}\colon V\to V_\mathbf{t}$ for the projection.
We will prove that $V_\mathbf{t}/\pi_\mathbf{t}(R)$ is finite of order $\prod_{\{p,q\}\in\mathbf{t}}|\Delta_{pq}|$.
It follows that $\pi_\mathbf{t}\colon R\to V_\mathbf{t}$ is injective, so that $\mathrm{torsion}(V/R)\xrightarrow{\pi_\mathbf{t}}V_\mathbf{t}/\pi_\mathbf{t}(R)$ is injective and consequently the order of $\mathrm{torsion}(V/R)$ divides $\prod_{\{p,q\}\in\mathbf{t}}|\Delta_{pq}|$.
Since by Lemma \ref{GCD_Lemma} the greatest common divisor of these quantities $\prod_{\{p,q\}\in\mathbf{t}}|\Delta_{pq}|$ is $1$, it follows that $\mathrm{torsion}(V/R)$ is the zero group, and the lemma is proved.

We now prove the claim that $V_\mathbf{t}/\pi_\mathbf{t}(R)$ is finite of order $\prod_{\{p,q\}\in\mathbf{t}}|\Delta_{pq}|$.
First, note that, if $\mathbf{s}$ is a tree on any finite set $A$ of size greater than $1$, then there is an $a\in A$ which lies in only one element $\{a,b\}$ of $\mathbf{s}$.
(This is a \emph{leaf} of the tree.)
Then $\mathbf{s}-\{\{a,b\}\}$ is a tree on $A-\{a\}$.
Repeatedly applying this process of deleting a leaf to the element $\mathbf{t}\in T_{[k+2]}$, we can write
\begin{gather*}
[k+2]=\{n_1,\ldots,n_{k+2}\},\\
\mathbf{t}=\{t_1,\ldots,t_{k+1}\},
\end{gather*}
where, for $p\leqslant(k+1)$, $n_p$ appears in $t_p$ but not in any of $t_{p+1},\ldots,t_{k+1}$.

Write $t_i=\{n_i,m_i\}$.
Then $V_{n_1m_1},\ldots,V_{n_{k+1}m_{k+1}}$ is a basis of $V_\mathbf{t}$, and $\pi_\mathbf{t} R_{n_1},\ldots,\pi_\mathbf{t} R_{n_{k+1}}$ span $\pi_\mathbf{t}(R)$.
Note that $\pi_\mathbf{t} V_{n_in}$ vanishes unless $\{n_i,n\}=t_j$ for some $j$, in which case we must have $j\leqslant i$.
If $j=i$ we must have $n=m_i$ and if $j<i$ we must have $n=n_j$, $n_i=m_j$.
So
\begin{eqnarray*}
\pi_\mathbf{t} R_{n_i}&=&\sum_{n\colon n\neq n_i}|\Delta_{n_in}|\pi_\mathbf{t} V_{n_in}\\
               &=&\sum_{n\colon\{n_i,n\}\in\mathbf{t}}|\Delta_{n_in}|V_{n_in}\\
               &=&|\Delta_{n_im_i}|V_{n_im_i}-\sum_{j<i\colon n_i=m_j}|\Delta_{m_jn_j}|V_{n_jm_j}.
\end{eqnarray*}
That is to say that, with respect to the chosen basis of $V_\mathbf{t}$, the matrix of relations is lower-triangular with diagonal entries $|\Delta_{n_im_i}|$.
Thus $V_\mathbf{t}/\pi_\mathbf{t}(R)$ has order $\prod|\Delta_{n_im_i}|=\prod_{\{p,q\}\in\mathbf{t}}|\Delta_{pq}|$, as required.
\end{proof}

\begin{lemma}\label{Determinant_Lemma}
Let $A$ be a finite set of size greater than $1$, ordered as $a_1,\ldots,a_n$.
Let $X_{ab}=X_{ba}$ be indeterminates indexed by unordered pairs of distinct elements $a,b\in A$.
Form the $(n-1)\times(n-1)$ matrix $\mathrm{M}$ with
\[\mathrm{M}_{lm}=\left\{\begin{array}{cl}\sum_{j\neq m+1}X_{a_ja_{m+1}}\ \ &\mathrm{if\ }l=m.\\ -X_{a_{l+1}a_{m+1}}\ \ &\mathrm{if\ }l\neq m.\end{array}\right.\]
That is,
\[\mathrm{M}=
\left(\begin{array}{cccc}
{\sum_{j\neq 2}X_{a_2a_j}}             & -X_{a_2a_3}                 & \cdots & -X_{a_2a_n} \\
-X_{a_2a_3}                            &  \sum_{j\neq 3}X_{a_3a_j}   &        &  \vdots         \\
\vdots                                 &                             & \ddots &  \vdots         \\
-X_{a_2a_n}                            &  \cdots                     & \cdots & {\sum_{j\neq n}X_{a_na_j}}
\end{array}\right)
.\]
Then
\[\det(\mathrm{M})=\sum_{\mathbf{t}\in T_A}\prod_{\{a,b\}\in\mathbf{t}}X_{ab}.\]
\end{lemma}

\begin{proof}
We will prove this lemma by induction on $|A|=n$, starting with $n=2$.
In the case $n=2$, the matrix $\mathrm{M}$ is just $(X_{a_1a_2})$, so that $\det(\mathrm{M})=X_{a_1a_2}$.
Since $T_A$ consists of the single tree $\{\{a_1,a_2\}\}$, the result
\[\det(\mathrm{M})=\sum_{\mathbf{t}\in T_A}\prod_{\{a,b\}\in\mathbf{t}}X_{ab}\]
is immediate.

Now assume that the claim is true whenever $A$ is replaced by a set of size less than $n$.

The determinant of $\mathrm{M}$ appears to depend on the ordering $a_1,\ldots,a_{n}$ of $A$, but this is not the case.
First, re-ordering $a_2,\ldots,a_n$ corresponds to conjugating $\mathrm{M}$ by a permutation matrix and so does not alter $\det(\mathrm{M})$.
Also, if we add the last $(n-2)$ rows of $\mathrm{M}$ to the first, and then the last $(n-2)$ columns of the resulting matrix to the first, and negate the first row and column, then we obtain the matrix analogous to $\mathrm{M}$ but for the ordering $a_2,a_1,a_3,\ldots,a_{n}$.
Thus $\det(\mathrm{M})$ is preserved under all permutations of $A-\{a_1\}$ and under the transposition of $a_1$ and $a_2$, and hence under all permutations of $A$.

It is clear that $\det(\mathrm{M})$ is a homogeneous polynomial of degree $(n-1)$.
Since each of the $X_{a_1a_j}$ appears in only one column of $\mathrm{M}$, no monomial in $\det(\mathrm{M})$ contains a square of such a generator; by symmetry, none of the generators $X_{ab}$ appears in any monomial of $\det(\mathrm{M})$ with multiplicity greater than one.
That is, $\det(\mathrm{M})$ is a sum of monomials of the form $\prod_{\{a,b\}\in\mathbf{t}}X_{ab}$ for $\mathbf{t}\in(A^{(2)})^{(n-1)}$.

Fix $\mathbf{t}\in (A^{(2)})^{(n-1)}-T_A$.
Since $\mathbf{t}$ is not a tree, the corresponding graph on $A$ is disconnected, and we may write $A=B\cup C$ with $B,C$ nonempty and disjoint, and $\mathbf{t}=\mathbf{u}\cup\mathbf{v}$ for $\mathbf{u}\subset B^{(2)}$ and $\mathbf{v}\subset C^{(2)}$.
Without loss of generality, $B=\{a_1,\ldots,a_r\}$ and $C=\{a_{r+1},\ldots,a_n\}$ for some $r$.
Setting $X_{a_ia_j}=0$ whenever $i\leqslant r$ and $j>r+1$, $\mathrm{M}$ becomes the partitioned matrix
$$\left(\begin{array}{c|c}
\mathrm{M}' & 0 \\
\hline 0 & \mathrm{M}''
\end{array}\right)$$
where $\mathrm{M}'$ is the analogue of $\mathrm{M}$ for the reduced set $\{a_1,\ldots,a_r\}$ and $\mathrm{M}''$ is a degenerate matrix with column-sum and row-sum zero.
Therefore, by setting $X_{ab}=0$ for $\{a,b\}\not\in\mathbf{t}$ (which as above involves setting $X_{a_ia_j}=0$ whenever $i\leqslant r$ and $j>r+1$), $\det(\mathrm{M})$ becomes 0.
It follows that $\det(\mathrm{M})$ does not contain a monomial of the form $\prod_{\{a,b\}\in\mathbf{t}}X_{ab}$.

The previous paragraph shows that $\det(\mathrm{M})$ contains no monomial of the form $\prod_{\{p,q\}\in\mathbf{t}}X_{pq}$ for $\mathbf{t}\in (A^{(2)})^{(n-1)}-T_A$, so that $\det(\mathrm{M})$ must be a sum of monomials indexed by $T_A$.
We will show that each monomial $\prod_{\{a,b\}\in\mathbf{t}}X_{ab}$ appears in $\det(\mathrm{M})$ with coefficient 1.
This will complete the proof.

Fix an element $\mathbf{t}\in T_A$.
Since $\mathbf{t}$ is a tree, there is an element of $A$ which appears in only one element of $\mathbf{t}$.
Without loss of generality, let this element be $a_n$, appearing only in $\{a_{n-1},a_n\}\in\mathbf{t}$.
Setting $X_{a_ia_n}=0$ for $i\leqslant (n-2)$, $\mathrm{M}$ becomes the partitioned matrix
$$\left(\begin{array}{c|c}
\mathrm{M}' & 0\\
\hline 0 & X_{a_{n-1}a_n}
\end{array}\right)$$
where $\mathrm{M}'$ is the analogue of $\mathrm{M}$ for the reduced set $A-\{a_n\}$.
By assumption, since $A-\{a_n\}$ has size less than $n$,
\[\det(\mathrm{M}')=\sum_{\mathbf{s}\in T_{A-\{a_n\}}}\prod_{\{a,b\}\in\mathbf{s}}X_{ab}.\]
If we set $X_{ab}=0$ for each $\{a,b\}\not\in\mathbf{t}$ (which involves setting $X_{a_i a_n}=0$ for $i\leqslant(n-2)$) then $\det(\mathrm{M}')$ becomes $\prod_{\{a,b\}\in\mathbf{t}-\{\{a_{n-1},a_n\}\}}X_{ab}$, and so $\det(\mathrm{M})$ becomes
\[X_{a_{n-1}a_n}\det(\mathrm{M}')=X_{a_{n-1}a_n}{\textstyle\prod_{\{a,b\}\in\mathbf{t}-\{\{a_{n-1},a_n\}\}}X_{ab}}=\prod_{\{a,b\}\in\mathbf{t}}X_{pq}.\]
That is, the monomial $\prod_{\{a,b\}\in\mathbf{t}}X_{ab}$ occurs in $\det(\mathrm{M})$ with coefficient 1.
The result follows.
\end{proof}

\begin{proof}[Proof of Proposition \ref{GOmega_Proposition}]
By Lemma \ref{NewDescription_Lemma}, and using the fact that $V_{pq}=\epsilon^p_qv_pv_q$, it follows that $G_\Omega=H^4(BT^k)/\mathrm{span}(\epsilon^p_qv_pv_q+\epsilon^q_rv_qv_r+\epsilon^r_pv_rv_p)$ is isomorphic to the group
\[V/(\mathrm{span}(R_p)+\mathrm{span}(V_{pq}+V_{qr}+V_{rp})).\]

Write $[k+2]=\{a_1,\ldots,a_{k+2}\}$.
Using this we choose the basis of $V$ given by $V_{a_ia_j}$ for $1\leqslant i<j\leqslant k+2$, ordered lexicographically in $i,j$.
Choose the basis $R_{a_2},\ldots,R_{a_{k+2}}$ of $\mathrm{span}(R_i)$, and the basis
$$V_{a_1a_i}+V_{a_ia_j}+V_{a_ja_1},\ \ 1\leqslant i<j\leqslant (k+2)$$
of $\mathrm{span}(V_{pq}+V_{qr}+V_{rp})$, ordered lexicographically in $i,j$.
Combining the second two we have an ordered basis of $\mathrm{span}(R_p)+\mathrm{span}(V_{pq}+V_{qr}+V_{rp})$, and we adjust it by replacing each $R_{a_p}$ with
$$R'_{a_p}=R_{a_p}-\sum_{q\neq p}|\Delta_{a_pa_q}|(V_{a_1a_p}+V_{a_pa_q}+V_{a_qa_1}).$$
Now, writing our basis of $\mathrm{span}(R_p)+\mathrm{span}(V_{pq}+V_{qr}+V_{rp})$ in terms of the basis of $V$, we obtain a \emph{matrix of relations} whose determinant has absolute value the order of the group $V/(\mathrm{span}(R_p)+\mathrm{span}(V_{pq}+V_{qr}+V_{rp}))=G_\Omega$.

It is easy to calculate the $R'_{a_p}$ and show that the matrix of relations has form
$$\left(\begin{array}{c|c}-M & 0\\
\hline ? & I\end{array}\right)$$
where
$$M=\left(\begin{array}{cccc}
{\sum_{j\neq 2}|\Delta_{a_2a_j}|}             & -|\Delta_{a_2a_3}|                 & \cdots & -|\Delta_{a_2a_{k+2}}| \\
-|\Delta_{a_2a_3}|                            &  \sum_{j\neq 3}|\Delta_{a_3a_j}|   &        &  \vdots             \\
\vdots                                        &                                    & \ddots &  \vdots             \\
-|\Delta_{a_2a_{k+2}}|                        &  \cdots                            & \cdots & {\sum_{j\neq k+2}|\Delta_{a_{k+2}a_j}|}
\end{array}\right).$$
The matrix $M$ is clearly the $\mathrm{M}$ of Lemma \ref{Determinant_Lemma}, with $A=[k+2]$ and $|\Delta_{pq}|$ in place of each $X_{pq}$, so that Lemma \ref{Determinant_Lemma} immediately gives us
\[\det(M)=\sum_{\mathbf{t}\in T_{[k+2]}}\prod_{\{p,q\}\in\mathbf{t}}|\Delta_{pq}|.\]
Consequently the matrix of relations has determinant
\[\pm \sum_{\mathbf{t}\in T_{[k+2]}}\prod_{\{p,q\}\in\mathbf{t}}|\Delta_{pq}|,\]
and the result follows.
\end{proof}

\section{Proof of Theorem \ref{Cohomology_Thm}.}\label{Cohomology_Section}
In this section we will prove Theorem \ref{Cohomology_Thm}.
To begin we prove a result on the cohomology of $N_\Omega$.
\begin{prop}\label{NOmegaCohomology_Prop}
$H^0(N_\Omega)=\Z$; $H^1(N_\Omega)$, $H^2(N_\Omega)$ vanish; $H^3(N_\Omega)$ is the free group on symbols $R_{pqr}$, for $p,q,r\in[k+2]$ distinct, subject only to antisymmetry in $p,q,r$ and
\[R_{pqr}+R_{qrs}+R_{rsp}+R_{spq}=0;\]
$H^4(N_\Omega)$ is a free group.
In the cohomology Serre spectral sequence $(E_r,d_r)$ for the bundle
\begin{equation}\label{NNOmegaBorel_Equation}N_\Omega\to\NN_\Omega\to BT^{k+2}\end{equation}
the classes $R_{pqr}$ are transgressive with $d_4(R_{pqr})=\rho_{pqr}$.
\end{prop}
\begin{proof}
The spectral sequence $(E_r,d_r)$ has $E_2^{\ast,\ast}=H^\ast(BT^{k+2})\otimes H^\ast(N_\Omega)$ and, by Theorem \ref{BorelCohomology_Thm},
\begin{align*}
E^{\ast,0}_\infty&=H^\ast(BT^{k+2})/\langle\rho_{pqr}\mid p,q,r\in[k+2]\ \mathrm{distinct}\rangle,\\
E^{\ast,j}_\infty&=0\ \ \mathrm{for\ }j>0.
\end{align*}
By Proposition \ref{NpOmega_Prop}, $N_\Omega$ is connected, so $H^0(N_\Omega)=\Z$.
Examining $(E_r,d_r)$ in total degree up to $4$, and using the description of $E_\infty$, we see that $H^1(N_\Omega)$ and $H^2(N_\Omega)$ do indeed vanish, and that there is a short exact sequence
\[0\to H^3(N_\Omega)\xrightarrow{d_4}H^4(BT^{k+2})\to H^4(BT^{k+2})/\langle\rho_{pqr}\rangle\to 0.\]
Consequently $d_4$ is an isomorphism of $H^3(N_\Omega)$ with the span of the $\rho_{pqr}$.
It is simple to show, using Lemma \ref{NewDescription_Lemma}, that the only relations among the $\rho_{pqr}$ are the analogues of the stated relations among the $R_{pqr}$, so the description of $H^3(N_\Omega)$ and $d_4$ follow.
Finally, there is an isomorphism
\[H^4(N_\Omega)\xrightarrow{d_2}\mathrm{ker}(d_4\colon H^2(BT^{k+2})\otimes H^3(N_\Omega)\to H^6(BT^{k+2})),\]
so that $H^4(N_\Omega)$ is free as required.
\end{proof}

\begin{proof}[Proof of Theorem \ref{Cohomology_Thm}]
The Borel construction of a group $G$ acting on a space $X$ is natural in $G$.
Therefore there is a commutative diagram whose rows are bundles and whose vertical maps are induced by $\Omega\colon T^k\to T^{k+2}$:
\[\xymatrix{
N_\Omega\ar[r]  &  \NN_\Omega\ar[r]   &  BT^{k+2}\\
N_\Omega\ar[r]\ar@{=}[u]  & ET^k\times_{T^k}N_\Omega\ar[r]\ar[u]^{E\Omega\times_{T^k}1} & BT^k\ar[u]^{B\Omega}
}\]
Moreover, $T^k$ acts freely on $N_\Omega$, so that $Q\colon ET^k\times_{T^k}N_\Omega\to T^k\backslash N_\Omega=S^7_\Omega$ is a homotopy-equivalence.
We can therefore use the cohomology Serre spectral sequence $(E'_r,d'_r)$ of the second row of the diagram above to compute $H^\ast(S^7_\Omega){\cong}H^\ast(ET^k\times_{T^k}N_\Omega)$; moreover, we can use the diagram to compare $(E'_r,d'_r)$ with $(E_r,d_r)$.

Note that ${E'}_2^{\ast,\ast}=H^\ast(BT^k)\otimes H^\ast(N_\Omega)$.
We display $E'_r$ in low degrees; a dot denotes the zero group, and $R=H^3(N_\Omega):$
\[\xymatrix @!=10 pt{
{\scriptstyle H^4(N_\Omega)}\ar[rrd]^{d'_4} &       &                    &       &             \\
R\ar[rrrrddd]^{d_4'}             & \cdot & {\scriptstyle H^2(BT^k)\otimes R} &       &             \\
\cdot         & \cdot & \cdot              & \cdot &             \\
\cdot         & \cdot & \cdot              & \cdot & \cdot       \\
\Z\ar@{-}'[r]'[rr]'[rrr]'[rrrr]\ar@{-}'[u]'[uu]'[uuu]'[uuuu]
   & \cdot & {\scriptstyle H^2(BT^k)}          & \cdot & {\scriptstyle H^4(BT^k)} & \cdot
}\]
Applying naturality of the Serre spectral sequence to the commutative diagram above, and using $B\Omega^\ast u_p=v_p$, we have
\[d'_4R_{pqr}=B\Omega^\ast\rho_{pqr}=\epsilon^p_qv_pv_q+\epsilon^q_rv_qv_r+\epsilon^r_pv_rv_p\]
so that $\mathrm{coker}(d'_4)=G_\Omega$.
Since $G_\Omega$ is finite this implies that $\rank (R)\geqslant\rank (H^4(BT^k))=\binom{k+1}{2}$; the relations among the $R_{pqr}$ described in Proposition \ref{NOmegaCohomology_Prop} show that $R$ is spanned by the $R_{1pq}$ for $1<p<q\leqslant (k+2)$, so that in fact $\rank (R)=\binom{k+1}{2}$.
It now follows, again since $G_\Omega$ is finite, that $\mathrm{ker}(d'_4)=0$.
Therefore ${E'}_\infty^{\ast,\ast}$ in low total degrees is
\[\xymatrix @!=10 pt{
{\scriptstyle\mathrm{ker}(d'_2)} &       &                    &       &             \\
\cdot              & \cdot &                    &       &             \\
\cdot              & \cdot & \cdot              &       &             \\
\cdot              & \cdot & \cdot              & \cdot &             \\
\Z\ar@{-}'[r]'[rr]'[rrr]'[rrrr]\ar@{-}'[u]'[uu]'[uuu]'[uuuu]
   & \cdot & {\scriptstyle H^2(BT^k)}          & \cdot & G_\Omega
}\]
Note that $\mathrm{ker}(d'_2)$ is a free group, so that $H^3(S^7_\Omega)=0$ and $H^4(S^7_\Omega)=\mathrm{ker}(d'_2)\oplus G_\Omega$.
However, by Poincar\'e duality $\rank H^4(S^7_\Omega)=\rank H^3(S^7_\Omega)$, so that in fact $\mathrm{ker}(d'_2)=0$ and we can conclude that in low degrees $H^\ast(S^7_\Omega)$ is given by:
\[\begin{array}{c|ccccc}
n & 0 & 1 & 2 & 3 & 4\\
H^n(S^7_\Omega) & \Z & 0 & \Z^k & 0 & G_\Omega
\end{array}\]
where the isomorphism of $H^2(S^7_\Omega)$ with $\Z^k$ is afforded by $c^\ast\colon H^2(BT^k)\to H^2(S^7_\Omega)$, where $c\colon S^7_\Omega\to BT^k$ is a classifying map for the action of $T^k$ on $N_\Omega$.
Poincar\'e duality now gives us the result.
\end{proof}

\section{The stable class of $TS^7_\Omega$.}\label{Pontrjagin_Section}
In this section we will prove Theorem \ref{Pontrjagin_Thm} by calculating the stable class of $TN_\Omega$ as an equivariant $T^{k+2}$ bundle, and then deducing the stable class of $TS^7_\Omega$ and its first Pontrjagin class.

\begin{proof}[Proof of Theorem \ref{Pontrjagin_Thm}]
Recall that $\mu_\Omega\colon S^{4k+7}\to\Im\H^k$ is $T^{k+2}$-invariant and has $\mathbf{0}$ as a regular value with $\mu^{-1}_\Omega(\mathbf{0})=N_\Omega$.
Therefore, if $\epsilon_\R$ denotes a trivial real line bundle (with trivial $T^{k+2}$ action), we have
\[TN_\Omega\oplus\epsilon_\R^{3k}\stackrel{T^{k+2}}{\cong} TS^{4k+7}|_{N_\Omega}.\]
(The symbol $\stackrel{G}{\cong}$ denotes $G$-equivariant isomorphism of bundles.)
Also
\[TS^{4k+7}\oplus\epsilon_\R\stackrel{T^{k+2}}{\cong}T\H^{k+2}|_{S^{4k+7}}\stackrel{T^{k+2}}{\cong}2\bigoplus_{p=1}^{k+2}L_p\]
where for $1\leqslant p\leqslant (k+2)$ $L_p$ denotes the trivial complex line bundle with $T^{k+2}$ action given by $\mathbf{t}\cdot z=t_pz$.
Therefore
$$TN_\Omega\oplus\epsilon_\R^{3k+1}\stackrel{T^{k+2}}{\cong}2\bigoplus_{p=1}^{k+2}L_p.$$
Now consider the action of $T^k$ on $N_\Omega$, which is given by the group homomorphism $\Omega\colon T^k\to T^{k+2}$.
We still have
\begin{equation}\label{Normal_Equation}TN_\Omega\oplus\epsilon_\R^{3k+1}\stackrel{T^k}{\cong}2\bigoplus_{p=1}^{k+2}L_p,\end{equation}
where now $L_p$ is the trivial complex line bundle with $T^k$-action given by $\mathbf{t}\cdot z=t_1^{\Omega_{p1}}\cdots t_k^{\Omega_{pk}}z$.
$T^k$ acts freely on $N_\Omega$ with quotient $S^7_\Omega$, so that taking the $T^k$ reduction of $TN_\Omega$ we obtain $T^k\backslash TN_\Omega\cong TS^7_\Omega\oplus\epsilon_\R^k$.
Also, for each $p$, $T^k \backslash L_p=l_p$, the complex line bundle over $S^7_\Omega$ with first Chern class $v_p$.
Thus, taking the $T^k$-reduction of \eqref{Normal_Equation},
\[TS^7_\Omega\oplus\epsilon_\R^{4k+1}\cong 2\bigoplus_{p=1}^{k+2}l_p\]
as required for the first part of the theorem.

Now $p_1(TS^7_\Omega)=-c_2(\mathbb{C}\otimes_\mathbb{R}TS^7_\Omega)=-c_2(\mathbb{C}\otimes_\mathbb{R}2{\textstyle\bigoplus}l_p),$ and
\begin{eqnarray*}
c(\mathbb{C}\otimes_\R 2{\textstyle\bigoplus}l_p)
&=& c({\textstyle\bigoplus}l_p)^2c({\textstyle\bigoplus}\bar l_p)^2\\
&=& \prod_{p=1}^{k+2}c(l_p)^2c(\bar l_p)^2\\
&=& 1-2\sum_{p=1}^{k+2}v_p^2,
\end{eqnarray*}
so that $p_1(S^7_\Omega)=2\sum_{p=1}^{k+2}v_p^2$, as required.
\end{proof}

\bibliographystyle{plain}
\bibliography{3spaper_bibliography}
\end{document}